\documentclass[12pt,a4paper]{article}

\title{Ordered Generating Systems of Finite Non-Abelian Groups }

\author{Robert Shwartz}

\begin{document}

\maketitle

\section{Introduction}

{\bf Definition:} Ordered Generating System: The elements $a_1,\cdots, 
a_n$ are considered Ordered Generating System of a group $G$, if every 
element  $g\in G$ has a unique representation in a form:

\vskip 0.3cm

$g=a_1^{i_1}a_2^{i_2}\cdots a_n^{i_n}$, where $0\leq i_k\leq m_k$, for 
some $m_k$, for every $1\leq k\leq n$. 

\vskip 0.3cm

We know from the basis theorem for finite abelian groups, that every 
abelian group has basis, and the basis by its definition, is the Ordered 
Generating System for an abelian group.

\vskip 0.2cm

Our motivation is generalizing the basis theorem, as it posible, for 
non-abelian groups. Hence we define the Ordered Generating System.  

\vskip 0.2cm

{\bf Lemma 1:} Let $G$ be a finite group, and let $H$ be it's normal 
subgroup. Assume $H$ has Ordered Generating system $a_1,\cdots, a_k$, and 
$G/H$ has Ordered Generating System $b_1H,\cdots b_lH$, then the elements 
$b_1\cdots b_l, a_1\cdots a_k$ are Ordered Generating System of $G$.

\vskip 0.2cm

{\bf Proof:} Since, every element of $G$ has a unique represntation in a 
form $ba$, where $bH\in G/H$, and $a\in H$, and since every element in 
$H$ has a unique representation in the form $a_1^{i_1}\cdots 
a_k^{i_k}$, and every element of $G/H$ has a unique representation in 
the form $b_1^{j_1}H\cdots b_l^{j_l}H$, we get that every element in $G$ 
has a unique representation in the form $b_1^{j_1}\cdots 
b_l^{j_l}a_1^{i_1}\cdots a_k^{i_k}$. Hence, the elements $b_1,\cdots, b_l, 
a_1, cdots, a_k$ are Ordered Generating System of the group $G$, by the 
definition of Ordered Generating System.

\vskip 0.3cm

{\bf Lemma 2:} Let $G$ be a finite group. Assume that each composition 
factor of $G$ has Ordered Generating System, then $G$ has Ordered 
Generating System, which is the union of the Ordered Generating Systems of 
the composition factors of $G$.

\vskip 0.2cm

{\bf Proof:} Let $G=G_0>G_1>G_2>\cdots>G_n=\{1\}$, be a sequence such that 
$G_{i+1}$ is a maximal normal subgroup of $G_i$, for every $1\leq i\leq 
n$. Then $G_i/G_{i+1}$ is isomorphic to one of the composition factor of 
$G$. Since $G_n=\{1\}$, $G_{n-1}$ is isomorphic to $G_{n-1}/G_n$ which is 
isomorphic to one of the composition factors of $G$. Hence, by the 
assumption of the Lemma, $G_{n-1}$ has Ordered Generating System. 
Since, $G_{n-2}/G_{n-1}$ is isomorphic to a composition factor of $G$, 
$G_{n-2}/G_{n-1}$ has Ordered Generating System, and since $G_{n-1}$ has 
Ordered Generating System, by Lemma 1, $G_{n-2}$ has Ordered 
Generating System.
Now assume by induction that every $G_i$ has Ordered Generating System, 
where $k\leq i\leq n-1$, and since $G_{k-1}/G_k$ has Ordered Generating 
System, by Lemma 1, $G_{k-1}$ has Ordered Generating System. Since 
$G=G_0$, by the induction $G$ has Ordered generating System.

\vskip 0.3cm

{\bf Known fact for solvable groups:} Since, the composition factor of a 
finite solvable group are cyclic groups, Lemma 2 implies that every 
solvable group has Ordered Genearating System. 
The Ordered Generating System is the elements which are coresponding to 
the generators of the cyclic groups in each composition factor. Hence, the 
existance of Ordered Generating System easily extendable to finite 
solvable groups.

\vskip 0.2cm

{\bf Non-Solvable groups:} Hence, we prove the existance of Ordered 
Generating System for some non-solvable groups.

\vskip 0.3cm

{\bf Lemma 3:} Let $G$ be a group. Assume $G$ has a subgroup $H$, such 
that $H$ has Ordered Generating System, and $gcd(|[G:H]|,|H|)=1$, and one 
of the following holds:

(i) $|[G:H]|=p^k$ where $p$ is prime number.

(ii) There exists an element of order $|[G:H]|$ in $G$.

(iii) There exists elements $a_1, a_2, \cdots a_n$, such that 
${a_i}^{m_i}\in H$, for $1\leq i\leq n$. Assume $m_1\cdot m_2\cdots 
m_n=|[G:H]|$, and all the elements of the form 
${a_1}^{i_1}{a_2}^{i_2}\cdots {a_n}^{i_n}\notin H$, where 
$0\leq i_k<m_k$.

Then $G$ has Ordered Generating System as well.

\vskip 0.2cm

{\bf Proof:} Assume (i) holds: Then the the $p$-sylow subgroup of $G$ does 
not belong to $H$, and we can take as a representative of the $p^k$ 
diferents cossets of $H$ in $G$, the $p^k$ diferent elements of the 
$p$-sylow subgroup of $G$. Since every $p$-group is a solvable group, by 
Lemma 2, the $p$-sylow subgroup of $G$ has Ordered Generating System.  
Then every element of $G$ has a unique representation in a form 
${a_1}^{i_1}{a_2}^{i_2}\cdots {a_m}^{i_m}{b_1}^{j_1}{b_2}^{j_2}\cdots 
{b_k}^{j_k}$, where $a_1,\cdots, a_m$ are the Ordered Generating System of 
the $p$-sylow subgroup of $G$, and $b_1,\cdots, b_k$ are the Ordered 
Generating System of the subgroup $H$ of $G$.

Assume (ii) holds: Then there exists an element $a$ of order $|[G:H]|$ in 
$G$. Since $gcd(|H|,|[G:H]|)=1$, $a^k\notin H$, for $1\leq k\leq |a|-1$. 
Hence, every element in $G$ has a unique representation of the form 
$a^ih$, where $h\in H$, and $0\leq i\leq |a|-1$. Then $a$, and the Ordered 
Generating System of $H$, is the Ordered Generating System of $G$.   

Assume (iii) holds: Then by the assumption of (iii) all the $|[G:H]|$ 
cossets of $H$ in $G$ can be written in the form 
${a_1}^{i_1}{a_2}^{i_2}\cdots {a_n}^{i_n}$, where $0\leq i_k<m_k$. Then, 
$a_1, a_2,\cdots, a_n$, and the Ordered Generating System of $H$ is the 
Ordered Generating System of $G$.

\vskip 0.3cm

We use the following Theorem:

\vskip 0.2cm

{\bf Theorem 1:} The groups, which composition factors are cyclic groups, 
$A_n$ or $PSL_n(q)$ has Ordered Generating system. 

\vskip 0.2cm

{\bf Proof:} By Lemma 2, it is enough to prove that every composition 
factor of $G$ has Ordered Generating System. Hence, the Theorem is 
interesting for the simple groups only. 

{\bf 1. The proof of the existance of Ordered Generating System for 
$A_n$:} The proof is by induction. $A_3$ is a cyclic group of order $3$, 
hence $A_3$ 
has Ordered Generating System. $A_4$ is a solvable group, hence by 
Lemma 2, $A_4$ has Ordered Generating System. Assume that $A_i$ has 
Ordered Generating System for every $i\leq 2k$.  $A_{2k+1}$ has a subgroup 
$H$ which is 
isomorphic to $A_{2k}$ (The stabilizer of $2k+1$). Let $A=(1,2,\cdots, 
2k+1)$. $A\in A_{2k+1}$, and every element of $A_{2k}$ has a unique 
representation of the form $A^{j}H$, where $0\leq j\leq 2k$. Hence, 
$A$, and the Ordered Generating System of $H$ (Which is exists by the 
assumption of the induction) is the Ordered Generating System of 
$A_{2k+1}$. Now, we prove that $A_{2k+2}$ has Ordered Generating System. 
Let $L$ be a subgroup of $A_{2k+2}$ which is the stabilizer 
of the point $2k+2$. Then $L$ is isomorphic to $A_{2k+1}$. 
Let $A=(1,2,\cdots, k, k+1)(k+2, k+3, \cdots, 2k+1, 2k+2)$, and let 
$B=(k+1, 2k+2)(1, 2k+1)$. Since the $2k+2$ elements $A^rB^s$, where $0\leq 
r\leq k$, $0\leq s\leq 1$ are taking the point $2k+2$ to the 2k+2 
diferent points in the permutation of 2k+2 points, then every element of 
$A_{2k+2}$ has a unique representation in a form $A^rB^sL$, where $0\leq 
r\leq k$, $0\leq s\leq 1$. Hence, $A, B$, and the Ordered Generating 
System of $L=A_{2k+1}$ are Ordered Generating System for $A_{2k+2}$. 
Hence, from the existance of Ordered Generating System for $A_{2k}$, we 
get that $A_{2k+1}$ and $A_{2k+2}$ have Ordered Generating System as well. 
Hence $A_n$ has Ordered Generating System for every $n$.
 
\vskip 0.2cm

{\bf 2. The proof of existance of Ordered Generating System for 
$PSL_n(q)$:} The proof is in induction in $n$. The order of $PSL_2(q)$ is 
$\frac{q(q-1)(q+1)}{2}$. $PSL_2(q)$ has a solvable subgroup $H$ of order 
$\frac{q(q-1)}{2}$, where $H$ is the subgroup corresponding to the upper 
triangular matrices. Since $H$ is solvable, By Lemma 2, $H$ has Ordered 
Generating System.  Since $H$ is corresponding to the upper 
triangular matrices, $H$ is the stabilizer of the point $\infty$ in 
the projective line over $F_q$. Since $|[PSL_2(q):H]|=q+1$, we apply Lemma 
3, for case (iii). Let $A$ be an element of order $\frac{q+1}{2}$ in 
$PSL_2(q)$, and let $B$ be an element of order $2$, such that $B$ is not 
corresponding to an upper triangular matrix in $PSL_2(q)$, and taking the 
point $\infty$ to a diferent point than every element $A^i$ (where $1\leq 
i<\frac{q+1}{2}$). Then, all the elements of the form $A^iB^j$ (where 
$0\leq i<\frac{q+1}{2}$, $0\leq j\leq 1$) are taking the point $\infty$ to 
the $q+1$ diferent points of the projective line $PF_q$. Hence there are 
$q+1$ diferent cosets of $H$ in $PSL_2(q)$ of the form $A^iB^j$, where 
$0\leq i<\frac{q+1}{2}$, $0\leq j\leq 1$. Then the elements $A$, $B$, and 
the Ordered Generating System of $H$ is Ordered Generating System for 
$PSL_2(q)$. 
 
Now assume that $PSL_{n-1}(q)$ has Ordered Generating System and we prove 
that $PSL_n(q)$ has Ordered Generating System as well. Let $H$ be a 
subgroup of $PSL_n(q)$ which is corresponding to the matrices where all 
the entries $a_{n,i}=0$, for $1\leq i\leq n-1$. Then the composition 
factors of $H$ are $PSL_{n-1}(q)$ and cyclic groups. $PSL_{n-1}(q)$ has 
Ordered Generating System by the assumption of the induction. Hence by 
Lemma 2, $H$ has Ordered Generating System. 
$|[PSL_n(q):H]|=\frac{q^n-1}{q-1}$. Since, $H$ is the subgroup which is 
the stabilizer of a subplane in the projective plane $PF_q^{n-1}$ and 
$PF_q^{n-1}$ contains $\frac{q^n-1}{q-1}$ points, we choose $A$ which  
order is $\frac{q^n-1}{(q-1)\cdot gcd(n,q-1)}$ and an element $B$ of order 
$gcd(n,q-1)$ in the case where $gcd(n,q-1)\neq 1$. 

\vskip 0.3cm

There are 5 sporadic Mathieu Groups: $M_{11}$, $M_{12}$, $M_{22}$, 
$M_{23}$, $M_{24}$.

\vskip 0.3cm

{\bf Theorem 3:} The Group $M_{11}$ has Ordered Generating System.

\vskip 0.2cm

{\bf Proof:} The sporadic group $M_{11}$ has order 7920, has a subgroup 
$H$ of index 11. The order of H is: 720. Since the composition factor of 
every non-solvable group of order $\leq 720$ is eithert $A_n$, $PSL_2(7)$, 
$PSL_2(8)$, or $PSL_2(11)$, by 2, $H$ has Ordered Generating System 
$a_2,\cdots a_n$. Since, the Order of $H$ is720. This order is prime to 
11, there is an element $a_1$ of Order 11 in $G$ which is not in $G$. 
Since, $[G:H]=11$, there are 11 diferent cosets of $H$ in $G$. Hence, 
$G=H\cup a_1H\cup \cdots \cup a_1^{10}H$. Then, $a_1$, and the Ordered 
Generating System $a_2, \cdots a_n$ of the subgroup $H$ are the Ordered 
Generating System of $G$.

\vskip 0.3cm

{\bf Theorem 4:} The Group $M_{12}$ has Ordered Generating System.

\vskip 0.2cm

{\bf Proof:} The Group $M_{12}$ is a subgroup of $S_{12}$ of Order 
$95040=2^6\cdot 3^3\cdot 5\cdot 11$, which is generared by:

$A=(1,2,3,4,5,6,7,8,9,10,11)$

$B=(5,6,4,10)(11,8,3,7)$

$C=(1,12)(2,11)(3,6)(4,8)(5,9)(7,10)$

The subgroup $H$ of $M_{12}$ which is generated by $A$ and $B$ is 
isomorphic to $M_{11}$, and then by Theorem 3, $H$ has Ordered Generating 
System.

Take the following elements: $X_1=A^9\cdot C\cdot 
A=(2,3,12)(1,8,4)(5,7,10)(6,9,11)$

$X_2=C=(1,12)(2,11)(3,6)(4,8)(5,9)(7,10)$

$X_3=A^8\cdot C\cdot A^3=(4,12)((3,5)((6,9)(7,11)(1,8)(2,10)$

Then the Ordered Generating System of $M_{12}$ is the Ordered Generating 
System of $H$ and the elements $X_1$, $X_2$, and $X_3$.

Since in $H$ is isomorphic to $M_{11}$, as a subgroup of $S_{12}$, where 
$H$ is the stabilizer of the point $12$  of the permutations in $S_{12}$. 
Then it can be shown easily that the 12 elements of the form 
$X_1^{i_1}\cdot X_2^{i_2}\cdot X_3^{i_3}$, 
where $0\leq i_1\leq 2$, $0\leq i_2\leq 1$, $0\leq i_3\leq 1$,  
are taking the point $12$ in $S_{12}$ to the 12 diferent points of 
$S_{12}$. Since $H$ is the stabilizer of the point $12$, every 
element in $M_{12}$ has a unique representation of the form $h\cdot 
X_1^{i_1}\cdot X_2^{i_2}\cdot X_3^{i_3}$, where $h\in H$, 
and $0\leq i_1\leq 2$, $0\leq i_2\leq 1$, $0\leq i_3\leq 1$. Since $H$ is
isomorphic to $M_{11}$, by Theorem 3, $H$ has Ordered Generating System.  
Then the Ordered Generating System of $M_{12}$ is: The Ordered Generating 
System of $H=M_{11}$, and the elements $X_1$, $X_2$, and $X_3$.

\vskip 0.3cm

{\bf Theorem 5:} The Group $M_{22}$ has Ordered Generating System.

\vskip 0.2cm

{\bf Proof:} $M_{22}$ is a group of Order $443520=2^7\cdot 3^2\cdot 5\cdot 
7\cdot 11$, is a subgroup of $S_{22}$ which is generated by the following 
3 permutations in $S_{22}$.

$X=(1,2,3,4,5,6,7,8,9,10,11)(12,13,14,15,16,17,18,19,20,21,22)$

$Y=(1,4,5,9,3)(2,8,10,7,6)(12,15,16,20,14)(13,19,21,18,17)$

$V=(11,22)(1,210(2,10,8,6)(12,14,16,20)(3,13,4,17)(5,19,9,18)$

\vskip 0.2cm

Let $H$ be a subgroup of $M_{22}$ which is the stabilizer of the point 
$22$ in the representation of $M_{22}$ as a subroup of $S_{22}$, which is 
generated by $X,Y,U$. Then $H$ is isomorphic to $PSL_3(4)$.

The 22 elements of the form $V^i\cdot X^j$, where $0\leq i\leq 1$, $0\leq 
j\leq 10$,  are taking the point $22$ in $S_{22}$ to the 22 diferent 
points of $S_{22}$. Since $H$ is the stabilizer of the point $22$, 
every element of $S_{22}$ has a unique representation of the form $h\cdot 
V^i\cdot X^j$, where $h\in H$,  $0\leq i\leq 1$, and $0\leq j\leq 10$. 
Since $H$ is isomorphic to $PSL_3(4)$, $H$ has Ordered Generating System 
by Theorem 2, and then the Ordered Generating System of $M_{22}$ are: The 
Ordered Generating System of $H$, and the elements $V$, and $X$.    

\vskip 0.3cm

{\bf Theorem 6:} The Group $M_{23}$ has Ordered Generating System.

\vskip 0.2cm

{\bf Proof:} $M_{23}$ is a group of Order $10200960=2^7\cdot 3^2\cdot 
5\cdot 7\cdot 11\cdot 23$. $M_{23}$ has a subgroup $H$ of index 23, which 
is isomorphic to $M_{22}$, and which order is prime to 23. By Theorem 5, 
$H$ has Ordered Generating System $a_2\cdot a_n$. Let $a_1$ be an element 
of order 23 in $M_{23}$. Since $[G:H]=23$, and the order of $H$ is prime 
to 23, there are 23 diferent cosets of $H$ in $G$ of the form $a_1^{i}H$, 
where $0\leq i<23$. Then, the elements $a_1$, and the Ordered Generating 
System of $H=M_{22}$, $a_2,\cdots a_n$,  are the Ordered Generating System 
of $M_{23}$. 

\vskip 0.3cm

{\bf Theorem 7:} The Group $M_{24}$ has Ordered Generating System.

\vskip 0.2cm

{\bf Proof:} $M_{24}$ is a group of Order $244823040=2^{10}\cdot 3^3\cdot 
5\cdot 7\cdot 11\cdot 23$. $M_{24}$ is a subgroup of $S_{24}$, which is 
generated by the following 3 permutations:

$D=(1,2,3,4,5,6,7,8,9,10,11,12,13,14,15,16,17,18,19,20,21,22,23)$

$E=(3,17,10,7,9)(4,13,14,19,5)(8,18,11,12,23)(15,20,22,21,16)$

$F=(1,24)(2,23)(3,12)(4,16)(5,18)(6,10)(7,20)(8,14)(9,21)(11,17)(13,22)(15,19)$

$M_{24}$ has a subgroup $H$ which is the stabilizer of the point $23$ in 
$S_{24}$, and isomorphic to $M_{23}$.

Let $X_1$ be $D^{-1}FD$, and let $X_2$ be $D^3F$. Then: 

$X_1=(2,24)(1,3)(4,13)(5,17)(6,19)(7,11)(8,21)(9,15)(10,22)(12,18)(14,23)(16,20)$
 
$X_2=(1,16,15,5,14,11,8,17,7,6,21,24)(2,18,9,3,10,22,23,12,19,13,4,20)$

Then the 24 diferent elements of the form $X_1^i\cdot X_2^j$, where $0\leq 
i<2$, and $0\leq j<12$,  are taking the point $24$ to the 24 diferent 
point of the permutations of $S_{24}$. Since $H$ is the stabilizer of the 
point $24$, every element in $S_{24}$ has a uniques representation in the 
form $h\cdot X_1^i\cdot X_2^j$, where $h\in H$, $0\leq i<2$, and $0\leq 
j<12$. Since $H$ is isomorphic to $M_{23}$, by Theorem 6, $H$ has Ordered 
Generating System. Then the Ordered Generating System of $H$ and the 
elements $X_1$, and $X_2$ are the Ordered Generating System of $M_{24}$. 

\end{document}